\documentclass[11pt,a4paper,leqno,twoside]{amsart}
\usepackage{latexsym,amssymb,amsmath}
\input xy
\xyoption{all}

\def\CC{\mathbb C}

\def\RR{\mathbb R}
\def\HH{\mathbb H}
\def\AA{{\mathbb A}}

\def\OO{\mathbb O}

\def\11{\mathbf 1}
\def\PP{\mathbb P}

\def\e1{\varepsilon_1}
\def\e2{\varepsilon_2}
\def\e3{\varepsilon_3}

\def\P2{{\PP}^2}

\def\00{\underline{0}}
\def\J0{{\cal J}_3(\underline{0})}

\def\PJ0{\PP({\cal J}_3(\underline{0}))}

\def\e{\varepsilon}

\def\AP2{{\AA\PP}^2}
\def\RP2{{\RR\PP}^2}
\def\CP2{{\CC\PP}^2}
\def\HP2{{\HH\PP}^2}
\def\OP2{{\OO\PP}^2}

\newtheorem{theo}{Theorem}[section]
\newtheorem{coro}[theo]{Corollary}
\newtheorem{lemm}[theo]{Lemma}
\newtheorem{prop}[theo]{Proposition}

\theoremstyle{definition}

\theoremstyle{remark}
\newtheorem{rema}[theo]{Remark}

\begin{document}
\title[Fields of dimension one]{On Fields of dimension one that are
Galois extensions of a global or local field}

\author{Ivan D. Chipchakov}
\address{Institute of Mathematics and Informatics\\Bulgarian Academy
of Sciences\\1113 Sofia, Bulgaria: E-mail address:
chipchak@math.bas.bg}

\begin{abstract}
Let $K$ be a global or local field, $E/K$ a Galois extension, and
Br$(E)$ the Brauer group of $E$. This paper shows that if $K$ is a
local field, $v$ is its natural discrete valuation, $v'$ is the
valuation of $E$ extending $v$, and $q$ is the characteristic of the
residue field $\widehat E$ of $(E, v')$, then Br$(E) = \{0\}$ if and
only if the following conditions hold: $\widehat E$ contains as a
subfield the maximal $p$-extension of $\widehat K$, for each prime
$p \neq q$; $\widehat E$ is an algebraically closed field in
case the value group $v'(E)$ is $q$-indivisible. When $K$ is a
global field, it characterizes the fields $E$ with Br$(E) = \{0\}$,
which lie in the class of tame abelian extensions of $K$.
We also give a criterion that, in the latter case, for any integer $n 
\ge 2$, there exists an $n$-variate $E$-form of degree $n$, which 
violates the Hasse principle.
\end{abstract}

\maketitle

\section{\bf Introduction}

\medskip
Let $F$ be a field, $F _{\rm sep}$ a separable closure of $F$, $F
_{\rm ab}$ the maximal abelian extension of $F$ in $F _{\rm sep}$,
$\mathbb{P}$ the set of prime numbers, and $F(p)$ the maximal
$p$-extension of $F$ in $F _{\rm sep}$, for each $p \in \mathbb{P}$.
We say that $F$ is a field of dimension $\le 1$, if the Brauer groups
Br$(F ^{\prime })$ are trivial, for all algebraic field extensions $F
^{\prime }/F$. It is known (cf. \cite{S1}, Ch.~II, 3.1) that dim$(F)
\le 1$ if and only if Br$(F ^{\prime }) = \{0\}$ when $F ^{\prime }$
runs across the set Fe$(F)$ of finite extensions of $F$ in $F _{\rm
sep}$. Note also that if dim$(F) \le 1$, then the absolute Galois
group $\mathcal{G}_{F} := \mathcal{G}(F _{\rm sep}/F)$ has
cohomological dimension cd$(\mathcal{G}_{F}) \le 1$ as a profinite
group; the converse holds in case $F$ is a perfect field.
\par
\medskip
It is well-known (cf. \cite{S1}, Ch. II, 3.1 and 3.2) that a field
$F$ satisfies dim$(F) \le 1$ whenever it is of type $C _{1}$ (or a
$C _{1}$-field), i.e. every $F$-form (a homogeneous nonzero
polynomial with coefficients in $F$) $f$ of degree deg$(f)$ in more
than deg$(f)$ variables has a nontrivial zero over $F$. The class of
$C _{1}$-fields contains finite fields (by Chevalley-Warning's
theorem, see, e.g., \cite{GiSz}, Theorem~6.2.6) as well as the
extensions of transcendency degree $1$ over any algebraically closed
field, by Tsen's theorem, and it is closed under taking algebraic
extensions (cf. \cite{L1}). Note also that if $F$ is a $C
_{1}$-field, then it is almost perfect, i.e. the following two
equivalent conditions hold: (i) every finite extension of $F$
possesses a primitive element; (ii) char$(F) = q \ge 0$, and in case
$q > 0$, the degree $[F\colon F ^{q}]$ of $F$ as an extension of its
subfield $F ^{q} = \{\beta ^{q}\colon \beta \in F\}$ is equal to $1$
or $q$. At the same time, almost perfect fields of dimension $\le 1$
need not be of type $C _{1}$. Indeed, the class $\mathcal{C}$
contains a field $E$ with dim$(E) \le 1$ that is not of type $C
_{1}$ in the following two cases: (i) $\mathcal{C}$ consists of
quasifinite fields of fixed characteristic $q$, where $q \in
\mathbb{P}$ or $q = 0$ (see \cite{Ax2}); (ii) $\mathcal{C}$ consists
of the algebraic extensions of an arbitrary global or local field
$K$; in addition, if $K$ is a global field and $v$ is its nontrivial
valuation, then $E$ can be chosen so as to have a Henselian
valuation extending $v$ \cite{Ch1}). The present research is
motivated by the following question.
\par
\medskip
{\bf Question 1.} Find whether a field $E$ with dim$(E) \le 1$ has
type $C _{1}$, provided that it is a Galois extension of a global or
local field $K$.
\par
\medskip
We pay particular attention to the case where the following
two conditions hold: (i) $E/K$ is an abelian extension, that is, the
Galois group $\mathcal{G}(E/K)$ is abelian; (ii) $E/K$ is tame, i.e.
for every finite extension $K ^{\prime }$ of $K$ in $E$ and any
nontrivial valuation $v$ of $K$, the characteristic char$(\widehat K
_{v})$ of the residue field $\widehat K _{v}$ of $(K, v)$ does not
divide the ramification index $e(K ^{\prime }/K) _{v}$ of $K ^{\prime
}/K$ relative to $v$. By definition, $e(K ^{\prime }/K) _{v}$ is the
index of the value group $v(K)$ as a subgroup of $v'(K ^{\prime })$,
for any valuation $v'$ of $K ^{\prime }$ extending $v$; since $E/K$
is abelian, whence, by Galois theory, so is $K ^{\prime }/K$, $e(K
^{\prime }/K) _{v}$ does not depend on the choice of $v'$ (see
\cite{L}, Ch. VI, Sect. 1, and Ch. XII, Corollary~6.3).
\par
\medskip
\section{\bf Statements of the main results}
\par
\medskip
For any field $K$ with a (nontrivial) Krull valuation $v$, $O
_{v}(K) = \{a \in K\colon \ v(a) \ge 0\}$ denotes the valuation ring
of $(K, v)$, $M _{v}(K) = \{\mu \in K\colon \ v(\mu ) > 0\}$ the
maximal ideal of $O _{v}(K)$, $O _{v}(K) ^{\ast } = \{u \in K\colon
\ v(u) = 0\}$ the multiplicative group of $O _{v}(K)$, $v(K)$ the
value group and $\widehat K _{v} = O _{v}(K)/M _{v}(K)$ the residue
field of $(K, v)$, respectively; $\overline {v(K)}$ is a divisible hull of
$v(K)$. We write for brevity $\widehat K$ instead of $\widehat K
_{v}$ when there is no danger of ambiguity. As usual, $v$ is said to be
discrete, if $v(K)$ is an infinite cyclic group. The valuation $v$ is
called Henselian, or else, we say that $(K, v)$ is a Henselian field,
if $v$ extends uniquely, up-to equivalence, to a valuation $v _{L}$
on each algebraic extension $L$ of $K$. The class of Henselian fields
contains every complete nontrivially real-valued field (see page
\pageref{k12}). When $v$ is Henselian, so is $v _{L}$, for any
algebraic field extension $L/K$. In this case, we put $O _{v}(L) = O
_{v _{L}}(L)$, $M _{v}(L) = M _{v_{L}}(L)$, $v(L) = v _{L}(L)$, and
denote by $\widehat L$ the residue field of $(L, v _{L})$; also, we
write $v$ instead of $v _{L}$ when the abbreviation is clear from the
context. As shown in \cite{Ch1}, if $(K, v)$ is a Henselian discrete
valued field (abbr., an HDV-field) with $\widehat K$ quasifinite, then
algebraic extensions of $K$ that are fields of dimension $\le 1$ are
characterized as follows:
\par
\medskip
(2.1) In order that dim$(L) \le 1$ it is necessary and sufficient
that the intersection $S(L) \cap \Sigma (L)$ be empty, where $S(L) =
\{p \in \mathbb{P}\colon \widehat L(p) \neq \widehat L\}$ and $\Sigma
(L)$ consists of those $p \in \mathbb{P}$, for which $v(L)$ is a
$p$-indivisible group.
\par
\medskip
The applicability of (2.1) is ensured by the former part of the
following result (cf. \cite{Ch1}, Lemma~4.1 and Corollary~4.3), which
indicates that the existence of fields of dimension $\le 1$ which are
not of type $C _{1}$ is not uncommon among algebraic extensions of
HDV-fields with finite residue fields:
\par
\medskip
(2.2) (a) For any HDV-field $(K, v)$ with $\widehat K$ quasifinite, and
any pair $S$, $\Sigma $ of subsets of $\mathbb{P}$, there exists an
algebraic extension $E/K$ satisfying $S(E) = S$ and $\Sigma (E) =
\Sigma $. The field $E$ has dimension $\le 1$ but is not of type $C
_{1}$, provided that $S \cap \Sigma = \emptyset $ and there is a pair
of integers $m _{j}$, $j = 1, 2$, such that $2 \le m _{1} \le m _{2}$,
$S$ contains all prime divisors of $m _{1}m _{2}$, and $\Sigma $
contains all prime divisors of $m _{1} + m _{2}$ (this requires that
$\gcd (m _{1}, m _{2}) = 1$). When dim$(E) \le 1$ and $E$ is not a $C
_{1}$-field, $S$ contains at least $2$ elements.
\par
(b) For each pair of integers $k _{1}$, $k _{2}$ with $2 \le k _{1}
< k _{2}$ and $\gcd (k _{1}, k _{2}) = 1$, there exist subsets $S
_{1}$ and $\Sigma _{1}$ of $\mathbb{P}$, such that $S _{1} \cap
\Sigma _{1} = \emptyset $, $S _{1}$ contains all prime divisors of $k
_{1}k _{2}$, and $\Sigma _{1}$ contains all prime divisors of $k _{1}
+ k _{2}$.
\par
\medskip
The main results of this paper are presented as two theorems. The
former one concerns the special case where $(K, v)$ is an HDV-field
with $\widehat K$ finite, and can be stated as follows:
\par
\medskip
\begin{theo}
\label{theo2.1} Let $(K, v)$ be an {\rm HDV}-field with $\widehat K$
finite and {\rm char}$(\widehat K) = q$, and suppose that $E$ is a
Galois extension of $K$, such that {\rm dim}$(E) \le 1$. Then
$\widehat E(p) = \widehat E$, for every $p \in \mathbb{P}$, $p \neq q$.
Moreover, if $v(E)$ is a $q$-indivisible group, then $\widehat E(q)
= \widehat E$; in particular, $\widehat E$ is algebraically closed.
\end{theo}
\par
\medskip
Theorem \ref{theo2.1} shows that, for every Galois extension $E/K$
satisfying dim$(E)$ $\le 1$, the set $S(E)$ is either empty
or equal to $\{q\}$. Note that, in the former case, it can be
deduced from Lang's theorem \cite{L1}, Theorem~10 (see also
Corollary \ref{coro4.2} below) that $E$ is of type $C _{1}$. On the
other hand, it seems that the sets $S(E ^{\prime })$ contain at
least $2$ elements, for all presently known algebraic extensions $E
^{\prime }$ of $K$ with dim$(E ^{\prime }) \le 1$, which are not $C
_{1}$-fields. Therefore, any approach to Question 1 should
rely on ideas and methods different from those used in \cite{Ch1}. It
should also be based on a satisfactory information on basic
algebraic and Diophantine properties of fields of dimension $\le 1$
that are Galois extensions of $K$ belonging to suitably chosen
special classes.
\par
\medskip
As a step in this direction, we consider fields of dimension $\le 1$
that are tame Galois extensions of a global field $K$ with abelian
Galois group. Frequently, we restrict to the special case where
char$(K) = 0$, i.e. $K$ is a number field. Our second main result is
contained in the following theorem.
\par
\medskip
\begin{theo}
\label{theo2.2} Let $K$ be a global field, $\Theta $ the compositum
of tame finite extensions of $K$ in $K _{\rm ab}$, and $E$ an
intermediate field of $\Theta /K$. Then:
\par
{\rm (a)} {\rm dim}$(E) \le 1$ if and only if $E$ is a nonreal field
and the residue fields of its discrete valuations are algebraically
closed; we have {\rm dim}$(\Theta ) \le 1$;
\par
{\rm (b)} If $K$ is a finite extension of the field $\mathbb{Q}$ of 
rational numbers, {\rm dim}$(E) \le 1$, $E$ contains as subfields 
$\Theta \cap K(2)$ and $\Theta \cap K(3)$, and $n \ge 2$ is an 
integer, then there exist $n$-variate $E$-forms $N _{m}$, $m \in 
\mathbb{N}$, of degree $n$, which are pairwise nonequivalent and 
violate the Hasse principle.
\end{theo}
\par
\medskip
The proof of Theorem \ref{theo2.2} (b) relies on the following two
facts concerning tame abelian extensions $E$ of an arbitrary global
field $K$ (see Corollary \ref{coro5.2} and Proposition \ref{prop5.3}
(b)):
\par
\medskip\noindent
(2.3) (a) Nontrivial Krull valuations of $E$ are discrete, and in
case dim$(E) \le 1$, their residue fields are algebraically closed;
\par
(b) If dim$(E) \le 1$, then for each $p \in \mathbb{P}$, there exist
infinitely many degree $p$ extensions $\Theta _{p,\nu }$, $\nu \in
\mathbb{N}$, of $K$ in $Y$.
\par
\medskip\noindent
Under the hypothesis of (2.3) (a), an algebraic extension $E
^{\prime }/E$ is said to be unramified, if $v'(E ^{\prime }) = v(E)$
whenever $v$ is a nontrivial valuation of $E$ and $v'$ is a
valuation of $E ^{\prime }$ extending $v$. When dim$(E) \le 1$, $E
^{\prime }/E$ is an unramified finite extension of degree $n$, $X
_{1}, \dots , X _{n}$ is a set of algebraically independent
variables over $E$, and $B$ is a basis of $E ^{\prime }$ as a vector
space over $E$, (2.3) implies the following:
\par
\medskip\noindent
(2.4) The norm form $N = N(X _{1}, \dots , X _{n})$ of $E ^{\prime
}/E$ associated with $B$ is an $n$-variate $E$-form of degree $n$,
which has a nontrivial zero over $E _{v}$, for each nontrivial
valuation $v$ of $E$, but does not possess a nontrivial $E$-zero. In
particular, $N$ violates the local-to-global principle over $E$.
\par
\medskip
Assuming that dim$(E) \le 1$ and $E$ is a tame extension of $K$ in
$K _{\rm ab}$, where $K/\mathbb{Q}$ is a finite Galois extension,
one formulates the main step towards the proof of Theorem
\ref{theo2.2} (b) as follows:
\par
\medskip
(2.5) For any integer $n \ge 2$, there exist unramified extensions
$E _{n,\nu }$, $\nu \in \mathbb{N}$, of $E$ in $E _{\rm ab}$, such 
that $[E _{n,\nu }\colon E] = n$ and $E _{n,\nu } \cap E _{n,\nu } 
^{\prime } = E$, for each index $\nu $, where $E _{n,\nu } ^{\prime 
}$ is the compositum of the fields $E _{n,\nu '}$, $\nu ' \neq \nu $; 
hence, by (2.3), every $E _{n,\nu }$ is embeddable as an 
$E$-subalgebra in $E _{w}$, for any nontrivial valuation $w$ of $E$. 
Moreover, if $N _{n,\nu }$ is a norm form of $E _{n,\nu }/E$, for 
each $\nu \in \mathbb{N}$, then $N _{n,\nu }$, $\nu \in \mathbb{N}$, 
are $n$-variate and pairwise nonequivalent $E$-forms of degree $n$, 
which violate the Hasse principle.
\par
\medskip
The former part of (2.5) indicates that, for each $\nu \in 
\mathbb{N}$, $N _{n,\nu }$ is a norm form of the field extension $E 
_{n,\nu }E _{n,\nu } ^{\prime }/E _{n,\nu } ^{\prime }$, whereas $N 
_{n,\nu '}$, $\nu ' \neq \nu $, decompose over $E _{n,\nu } ^{\prime 
}$ into products of $n$ linear forms in $n$ variables; in particular, 
this implies the latter part of (2.5).
\par
\medskip
The basic notation, terminology and conventions kept in this paper
are standard and virtually the same as in \cite{TW}, \cite{L} and
\cite{CF}. Throughout, Brauer and value groups are written
additively, Galois groups are viewed as profinite with respect to
the Krull topology, and by a profinite group homomorphism, we mean a
continuous one. An additively written abelian group $A$ is called
reduced, if it does not possess a nonzero divisible subgroup. When
$A$ is also a torsion group, $A _{p}$ denotes its $p$-component, for
every $p \in \mathbb{P}$. For any field extension $E ^{\prime }/E$,
we write I$(E ^{\prime }/E)$ for the set of intermediate fields of $E
^{\prime }/E$, $\pi _{E/E'}$ for the scalar extension map Br$(E) \to
{\rm Br}(E ^{\prime })$, and Br$(E ^{\prime }/E)$ for the relative
Brauer group of $E ^{\prime }/E$ (the kernel of $\pi _{E/E'}$). When
$E ^{\prime }/E$ is a Galois extension, $\mathcal{G}(E ^{\prime
}/E)$ denotes its Galois group and $X(E ^{\prime }/E)$ is the
continuous character group of $\mathcal{G}(E ^{\prime }/E)$. As
usual, $\mathbb{Z} (p ^{\infty })$ stands for a quasi-cyclic
$p$-group, $\mathbb{Z} _{p}$ is the ring of $p$-adic integers, and a
$\mathbb{Z} _{p}$-extension means a Galois extension $\Psi ^{\prime
}/\Psi $ with $\mathcal{G}(\Psi ^{\prime }/\Psi )$ isomorphic to the
additive group of $\mathbb{Z} _{p}$. For any field $E$, $E ^{\ast }$
is its multiplicative group, and $E ^{\ast n} = \{a ^{n}\colon a \in E
^{\ast }\}$, for each $n \in \mathbb{N}$. The field $E$ is said to be
formally real, if $-1$ is not presentable as a finite sum of elements
of $E ^{\ast 2}$; $E$ is called a nonreal field, otherwise. The value
group of any discrete valued field is assumed to be an ordered
subgroup of the additive group of the field $\mathbb{Q}$; this is
done without loss of generality, in view of \cite{E3},
Theorem~15.3.5, and the fact that $\mathbb{Q}$ is a divisible hull of
its infinite subgroups (see page \pageref{k99}).
\par
\medskip
Here is an overview of this paper: Section 3 includes
valuation-theoretic preliminaries used in the sequel. Theorem
\ref{theo2.1} is proved in Section 4. It is also shown there that if
$E$ is a field with dim$(E) \le 1$, which is an abelian extension of
a local field $K$, then the valuation of $E$ extending the natural
valuation of $K$ is discrete, and its residue field $\widehat E$ is
algebraically closed. Theorem \ref{theo2.2} (a) and statements (2.3)
are proved in Section 5, where we consider fields $\Delta $ of
dimension $\le 1$ that are abelian tame extensions of a global field
$K$. We prove (2.3) (b) by showing that, in characteristic zero, the
groups $X(\Delta /K)$ are reduced with finitely many elements of
infinite height, and with infinitely many elements of order $p$, for
each $p \in \mathbb{P}$. The latter property of $X(\Delta /K)$ is
used in Section 6 for proving Theorem \ref{theo2.2} (b). In Section
5, we also characterize the maximal abelian tame extension of
$\mathbb{Q}$ in $\mathbb{Q} _{\rm sep}$ as the extension of
$\mathbb{Q}$ obtained by adjunction of the primitive $p$-th root of
unity $\varepsilon _{p} \in \mathbb{Q} _{\rm sep}$, for all $p \in
\mathbb{P}$; its analog in characteristic $q > 0$ is the rational
function field $\mathbb{F} _{q,{\rm sep}}(X)$.
\par
\medskip
\section{\bf Preliminaries and characterizations of algebraic
extensions $E$ of local or global fields with Br$(E) _{p} \neq \{0\}$, for a
given prime $p$}
\par
\medskip
\label{k12} Let $K$ be a field with a (nontrivial) Krull valuation
$v$. It is known (cf. \cite{L}, Ch. XII) that $v$ is Henselian in
case $v(K)$ embeds as an ordered subgroup in the additive group
$\mathbb R$ of real numbers and $K$ has no proper separable
extensions in its completion $K _{v}$ (with respect to the topology
of $K$ induced by $v$). For an arbitrary $v$, the Henselian condition
has the following two equivalent forms (cf. \cite{E3}, Sect. 18.1,
and \cite{L}, Ch. XII, Sect. 4):
\par
\medskip\noindent
(3.1) (a) Given a polynomial $f(X) \in O _{v}(K) [X]$ and an element
$a \in O _{v}(K)$, such that $2v(f ^{\prime }(a)) < v(f(a))$, where
$f ^{\prime }$ is the formal derivative of $f$, there is a zero $c
\in O _{v}(K)$ of $f$ satisfying the equality $v(c - a) = v(f(a)/f
^{\prime }(a))$;
\par
(b) For each normal extension $\Omega /K$, $v ^{\prime }(\tau (\mu ))
= v ^{\prime }(\mu )$ whenever  $\mu \in \Omega $, $v ^{\prime }$ is
a valuation of $\Omega $ extending $v$, and $\tau $ is a
$K$-automorphism of $\Omega $.
\par
\medskip
When $v$ is Henselian, so is $v _{L}$, for any algebraic field
extension $L/K$. In this case, we put $O _{v}(L) = O _{v _{L}}(L)$,
$M _{v}(L) = M _{v_{L}}(L)$, $v(L) = v _{L}(L)$, and denote by
$\widehat L$ the residue field of $(L, v _{L})$; also, we write $v$
instead of $v _{L}$ when there is no danger of ambiguity. Clearly,
$\widehat L/\widehat K$ is an algebraic extension and $v(K)$ is an
ordered subgroup of $v(L)$, such that $v(L)/v(K)$ is a torsion
group; hence, one may assume without loss of generality that
\label{k99} $v(L)$ is an ordered subgroup of $\overline {v(K)}$. By
Ostrowski's theorem (cf. \cite{E3}, Theorem~17.2.1), if $[L\colon
K]$ is finite, then it is divisible by $[\widehat L\colon \widehat
K]e(L/K)$ and $[L\colon K][\widehat L\colon \widehat K] ^{-1}e(L/K)
^{-1}$ has no divisor $p \in \mathbb{P}$, $p \neq {\rm
char}(\widehat K)$; here $e(L/K)$ is the index of $v(K)$ in $v(L)$.
Ostrowski's theorem implies the following:
\par
\medskip
(3.2) The quotient groups $v(K)/pv(K)$ and $v(L)/pv(L)$ are
isomorphic, if $p \in \mathbb{P}$ and $[L\colon K] < \infty $. When
char$(\widehat K) \nmid [L\colon K]$, the natural embedding of $K$
into $L$ induces canonically an isomorphism $v(K)/pv(K) \cong
v(L)/pv(L)$.
\par
\medskip
The extension $L/K$ is defectless, i.e. $[L\colon K] =
[\widehat L\colon \widehat K]e(L/K)$, in the following three cases:
\par
\medskip\noindent
(3.3) (a) If char$(\widehat K) \nmid [L\colon K]$ (apply Ostrowski's
theorem).
\par
(b) If $(K, v)$ is HDV and $L/K$ is separable (see \cite{E3}, Sect.
17.4).
\par
(c) If $(K, v)$ is HDV and $K$ is almost perfect. In particular,
this holds in the following two cases: when $(K, v)$ is a complete
discrete valued field with $\widehat K$ perfect; if $(K, v)$ is HDV
and $K$ is an algebraic extension of a global field.
\par
\medskip\noindent
Assume that $(K, v)$ is a nontrivially valued field and $R$ is a
finite extension of $K$, which has a unique, up-to equivalence,
valuation $v _{R}$ extending $v$. The extension $R/K$ is called
inertial relative to $v$, if the residue field $\widehat R$ of $(R, v
_{R})$ is separable over $\widehat K$, and $[R\colon K] = [\widehat
R\colon \widehat K]$; we say that $R/K$ is totally ramified with
respect to $v$, if the index of $v(K)$ in $v _{R}(R)$ equals
$[R\colon K]$. When $v$ is Henselian, $R/K$ is totally ramified, if
$e(R/K) = [R\colon K]$. Under the same condition, inertial extensions
of $K$ (with respect to $v$) have the following useful properties
(see \cite{TW}, Theorem~A.23):
\par
\medskip\noindent
(3.4) (a) An inertial extension $R ^{\prime }/K$ is Galois if and
only if so is $\widehat R ^{\prime }/\widehat K$. When this holds,
$\mathcal{G}(R ^{\prime }/K)$ and $\mathcal{G}(\widehat R ^{\prime
}/\widehat K)$ are canonically isomorphic.
\par
(b) The compositum $K _{\rm ur}$ of inertial extensions of $K$ in $K
_{\rm sep}$ is a Galois extension of $K$ with $v(K _{\rm ur}) =
v(K)$ and $\mathcal{G}(K _{\rm ur}/K) \cong \mathcal{G}_{\widehat
K}$.
\par
(c) Finite extensions of $K$ in $K _{\rm ur}$ are inertial, and the
natural mapping of $I(K _{\rm ur}/K)$ into $I(\widehat K _{\rm
sep}/\widehat K)$, by the rule $L \to \widehat L$, is bijective.
\par
\medskip\noindent
Returning to the case of $(K, v)$ Henselian, we denote by $K _{\rm
tr}$ the compositum of those $\Lambda \in {\rm Fe}(K)$, which are
tamely ramified extensions of $K$, i.e. $e(\Lambda /K)$ is not
divisible by char$(\widehat K)$. It is known that $K _{\rm tr}/K$ is
a Galois extension with the following properties (see \cite{TW},
Appendix~A2):
\par
\medskip
(3.5) (a) $K _{\rm ur} \in I(K _{\rm tr}/K)$ and $K _{\rm tr}/K
_{\rm ur}$ is Galois with $\mathcal{G}(K _{\rm tr}/K _{\rm ur})$
isomorphic to the topological group product $\prod _{p \in \Sigma }
\mathbb{Z} _{p} ^{\tau (p)}$, where $\Sigma $ is the set $\{p \in
\mathbb{P}\colon v(K) \neq pv(K)\}$ and for each $p \in \Sigma $,
$\mathbb{Z} _{p} ^{\tau (p)}$ is the topological product of
isomorphic copies of $\mathbb{Z} _{p}$, indexed by a set of
cardinality equal to the dimension $\tau (p)$ of the group
$v(K)/pv(K)$, viewed as a vector space over $\mathbb{F}_{p}$ (with
respect to the operations naturally induced by the addition in
$v(K)$);
\par
(b) If $K _{\rm sep} \neq K _{\rm tr}$, then char$(\widehat K) = q
> 0$ and $K _{\rm sep}/K _{\rm tr}$ is a $q$-extension.
\par
\medskip
\begin{lemm}
\label{lemm3.1} Let $(K, v)$ be a real-valued field, $(K _{v}, \bar
v)$ its completion, and $(K ^{\prime }, v')$ an intermediate valued
field of $(K _{v}, \bar v)/(K, v)$. Assume that $(K ^{\prime }, v')$
is Henselian and identify $K ^{\prime } _{\rm sep}$ with its $K
^{\prime }$-isomorphic copy in $K _{v,{\rm sep}}$. Then:
\par
{\rm (a)} $K ^{\prime } _{\rm sep} \cap K _{v} = K ^{\prime }$, and
each $\Lambda \in {\rm Fe}(K _{v})$ contains a primitive element
$\lambda \in K ^{\prime } _{\rm sep}$ over $K _{v}$, such that $[K
_{v}(\lambda )\colon K _{v}] = [K ^{\prime }(\lambda )\colon K
^{\prime }]$;
\par
{\rm (b)} $K ^{\prime } _{\rm sep}K _{v} = K _{v,{\rm sep}}$ and
$\mathcal{G}_{K'} \cong \mathcal{G}_{K _{v}}$;
\par
{\rm (c)} The mapping, say $f$, of {\rm Fe}$(K ^{\prime })$ into
{\rm Fe}$(K _{v})$, by the rule $\Lambda ^{\prime } \to \Lambda
^{\prime }K _{v}$, is bijective and degree-preserving. Moreover, $f$
and the inverse mapping $f ^{-1}: {\rm Fe}(K _{v}) \to {\rm Fe}(K
^{\prime })$, preserve the Galois property and the isomorphism class
of the corresponding Galois groups;
\par
{\rm (d)} For each $\nu \in \mathbb{N}$ not divisible by {\rm
char}$(K)$, $K _{v} ^{\ast } = K ^{\prime \ast }K _{v} ^{\ast \nu }$.
\end{lemm}
\par
\medskip
\begin{lemm}
\label{lemm3.2} In the setting of Lemma \ref{lemm3.1}, suppose
that $(K, v)$ is Henselian, identify $K _{\rm sep}$ with its
$K$-isomorphic copy in $K _{v,{\rm sep}}$, fix an extension $R$ of
$K$ in $K _{\rm sep}$, and put $R ^{\prime } = K ^{\prime }R$. Then
$(R ^{\prime }, v'_{R'})$ is an intermediate valued field of $(R
_{v}, \bar v_{R})/(R, v_{R})$.
\end{lemm}
\par
\medskip
Our next lemma characterizes those fields of dimension $\le 1$, which
lie in the class of algebraic extensions of any HDV-field $(K, v)$
with $\widehat K$ qusifinite. The lemma has been proved in \cite{Ch1}
(see \cite{S1}, Ch. II, for a proof in case $(K, v)$ is a local field
with char$(K) = 0$). In particular, it shows that if $(K, v)$ is HDV,
$\widehat K$ is quasifinite, and $E/K$ is an algebraic extension,
then dim$(E) \le 1$ if and only if the intersection $S(E) \cap \Sigma
(E)$ is empty, where $S(E) = \{p \in \mathbb{P}\colon
\mathcal{G}(\widehat E(p)/\widehat E) \cong \mathbb{Z} _{p}\}$ and
$\Sigma (E) = \{p \in \mathbb{P}\colon v(E) \neq pv(E)\}$.
\par
\medskip
\begin{lemm}
\label{lemm3.3} Assume that $(K, v)$ is an {\rm HDV}-field with
$\widehat K$ quasifinite, fix some $p \in \mathbb{P}$, and take an
algebraic field extension $R/K$. Then {\rm Br}$(R) _{p} = \{0\}$ if
and only if each of the following three equivalent conditions is fulfilled:
\par
{\rm (a)} {\rm Br}$(R ^{\prime }) _{p} = \{0\}$, for every algebraic
extension $R ^{\prime }/R$; this holds if and only if Br$(R ^{\prime
}) _{p} = \{0\}$, when $R ^{\prime }$ runs across the set Fe$(R)$;
\par
{\rm (b)} $p$ does not divide the period of the quotient group $R
^{\prime \ast }/N(R _{1} ^{\prime }/R ^{\prime })$, for any $R _{1}
^{\prime } \in {\rm Fe}(R)$ and every $R ^{\prime } \in I(R _{1}
^{\prime }/R)$;
\par
{\rm (c)} $v(R) = pv(R)$ or $\widehat R(p) = \widehat R$.
\par\noindent
In addition, if $R/K$ is separable or $K$ is almost perfect, or $p
\neq {\rm char}$(K), then {\rm Br}$(R) _{p} = \{0\}$ if and only if
there are finite extensions $R _{n}$, $n \in \mathbb{N}$, of $K$ in
$R$, such that $p ^{n}$ divides $[R _{n}\colon K]$, for each index
$n$.
\end{lemm}
\par
\medskip
The following lemma characterizes fields of dimension $\le 1$ that
are algebraic extensions of a global field as follows:
\par
\medskip
\begin{lemm}
\label{lemm3.} Let $K$ be a global field, $E/K$ an algebraic
extension, $\mathcal{M}_{E}$ the set of equivalence classes of
nontrivial Krull valuations of $E$, and $R _{E}$ a system of
representatives of $\mathcal{M}_{E}$ consisting of
$\mathbb{Q}$-valued valuations. For each $v' \in R _{E}$, fix a
completion $E _{v'}$ of $E$ with respect to the topology of $v'$,
denote by $v$ the valuation of $K$ induced by $v'$, and identifying
$K _{v}$ with the closure of $K$ in $E _{v'}$, put $E(v') = E.K
_{v}$. Then {\rm Br}$(E) _{p} = \{0\}$, for some $p \in \mathbb{P}$,
if and only if {\rm Br}$(E(v')) _{p} = \{0\}$, $v' \in R _{E}$, and
in case $p = 2$, $E$ is a nonreal field.
\end{lemm}
\par
\medskip
Lemma \ref{lemm3.} is implied by Lemma \ref{lemm3.3} and
\cite{S1}, Ch. II, Proposition~9.
\par
\medskip
\begin{lemm}
\label{lemm3.4} Assume that $K$ is a global field, $\mathcal{M}_{K}$
is the set of equivalence classes of discrete valuations of $K$, and
$R _{K}$ is a system of representatives of $\mathcal{M}_{K}$. Let
$E$ be a Galois extension of $K$, and for each $v \in R _{K}$, let
$E(v)$ be a $K$-isomorphic copy of $E$ in $K _{v,{\rm sep}}$. Then
{\rm Br}$(E) _{p} = \{0\}$, for some $p \in \mathbb{P}$, if and only
if $E$ is a nonreal field with {\rm Br}$(E(v)K _{v}) _{p} = \{0\}$,
for every $v \in R _{K}$.
\end{lemm}
\par
\medskip
For any field $R$, it is known that conditions (a) and (b) of Lemma
\ref{lemm3.3} are equivalent, and if they hold, then
$\mathcal{G}_{R}$ is a profinite group of cohomological
$p$-dimension cd$_{p}(\mathcal{G}_{R}) \le 1$; this implication is
an equivalence in case $R$ is perfect or $p \neq {\rm char}(R)$ (cf.
\cite{S1}, Ch. II, 3.1, and \cite{GiSz}, Theorem~6.1.8). Generally,
the condition that Br$(R) _{p} = \{0\}$ is weaker than conditions
(a) and (b) of Lemma \ref{lemm3.3}. It is known, however, that if $E
\in I(\mathbb{Q} _{\rm sep}/\mathbb{Q})$ and Br$(E) _{p} = \{0\}$,
for some $p \in \mathbb P$, then Br$(E ^{\prime }) _{p} = \{0\}$
whenever $E ^{\prime } \in I(\mathbb{Q} _{\rm sep}/E)$ \cite{FS},
Theorem~4; this also follows from Lemma \ref{lemm3.3} and \cite{S1},
Ch. II, Proposition~9).
\par
\medskip
\section{\bf Proof of Theorem \ref{theo2.1} and applications to
fields of dimension $\le 1$ that are tamely ramified abelian
extensions of a local field}
\par
\medskip
Our first objective in this Section is to prove Theorem
\ref{theo2.1}. One may assume without loss of generality that $E \in
I(K _{\rm sep}/K)$. As in Section 4, we put $S(L) = \{p \in
\mathbb{P}\colon \widehat L(p) \neq \widehat L\}$ and $\Sigma (L) =
\{p \in \mathbb{P}\colon v(L) \neq pv(L)\}$, for every $L \in I(K
_{\rm sep}/K)$. Since dim$(E) \le 1$, we have $S(E) \cap \Sigma (E)
= \emptyset $, so it suffices to prove that $\widehat E(p) =
\widehat E$, for an arbitrary fixed $p \in \Sigma (E)$, $p \neq {\rm
char}(\widehat K)$. Fix a primitive $p$-th root of unity $\varepsilon
\in K _{\rm sep}$. It is well-known that $K(\varepsilon )/K$ is a
Galois extension, $\mathcal{G}(K(\varepsilon )/K)$ is cyclic and
$[K(\varepsilon )\colon K] \mid p - 1$ (cf. \cite{L}, Ch. VI, Sect.
3). This implies $E(\varepsilon )/K$ is Galois and $[E(\varepsilon
)\colon E] \mid [K(\varepsilon )\colon K]$, so it follows from (3.2)
that $\Sigma (E(\varepsilon )) = \Sigma (E)$. Observing also that $L
_{\rm ur} = K _{\rm ur}L$ and $\mathcal{G}_{\widehat K} \cong \prod
_{\ell \in \mathbb{P}} \mathbb{Z} _{\ell }$, for each $L \in I(K
_{\rm sep}/K)$, one obtains from (3.4) (b) and Galois theory that
$\widehat L(\ell )/\widehat L$ is a $\mathbb{Z} _{\ell }$-extension,
for every $\ell \in S(L)$. It is therefore clear that
$S(E(\varepsilon )) = S(E)$ and it suffices to prove that $p \notin
S(E)$ in the special case where $\varepsilon \in E$. Consider now
the fields $E _{0} = E \cap K _{\rm ur}$ and $E _{0} ^{\prime } = E
\cap K _{\rm tr}$, and fix a system $\varepsilon _{n}$, $n \in
\mathbb{N}$, of elements of $K _{\rm sep}$, chosen so that
$\varepsilon _{1} ^{p} = \varepsilon $ and $\varepsilon _{n} ^{p} =
\varepsilon _{n-1}$, for every index $n \ge 2$. It is easily
verified that $\varepsilon _{n} \notin K$ whenever $n$ is
sufficiently large, which ensures that the extension $K ^{\prime }$
of $K$ generated by the set $\{\varepsilon _{n}\colon n \in
\mathbb{N}\}$ is a $\mathbb{Z} _{p}$-extension of $K$ and equals the
field $K(p) \cap K _{\rm ur}$. Thus it turns out that the equality
$\widehat E(p) = \widehat E$ will follow, if we show that $K
^{\prime }$ is included in $E$, i.e. $\varepsilon _{n} \in E$, for
all $n \in \mathbb{N}$. Note also that $p \in s(E)$ if and only if
$p \in S(E _{0} ^{\prime })$. This can be deduced from (3.2) and the
well-known fact that $\mathcal{G}(E/E _{0} ^{\prime })$ is a
pro-$q$-subgroup of $\mathcal{G}(E/K)$, where $q = {\rm
char}(\widehat K)$. In other words, one may assume for the rest of
our proof that $E = E _{0} ^{\prime }$, i.e. $E \in I(K _{\rm
tr}/K)$.
\par
Our argument also relies on the following two facts: (i) finite
extensions of $E _{0}$ in $E$ are tamely and totally ramified; (ii)
$v(E _{0}) = v(K)$, by (3.4) (b), i.e. $(E _{0}, v)$ is an
HDV-field. Note further that $\mathcal{G}(E/E _{0})$ is an abelian
group. Indeed, it follows from Galois theory and the equality $E
_{0} = E \cap K _{\rm ur}$ that $E.K _{\rm ur}/K _{\rm ur}$ is a
Galois extension with $\mathcal{G}(E.K _{\rm ur}/K _{\rm ur}) \cong
\mathcal{G}(E/E _{0})$. Since $E.K _{\rm ur} \in I(K _{\rm tr}/K
_{\rm ur})$, one also sees that $\mathcal{G}(E/E _{0})$ is a
homomorphic image of $\mathcal{G}(K _{\rm tr}/K _{\rm ur})$, so it
follows from (3.5) (a) that $\mathcal{G}(E/E _{0})$ is abelian.
Furthermore, (3.5) (a) and the cyclicity of the group $v(E _{0}) =
v(K)$ imply $\mathcal{G}(E/E _{0})$ is a procyclic group. Hence, by
(3.2) and the assumption that $v(E) = pv(E)$, $E$ contains as a
subfield a $\mathbb{Z} _{p}$-extension $\Lambda _{\infty }$ of $E
_{0}$. This means that $E _{0}$ has (Galois) extensions $\Lambda
_{n}$, $n \in \mathbb{N}$, such that $\Lambda _{\infty } = \cup
_{n=1} ^{\infty } \Lambda _{n}$, and for each index $n$, $\Lambda
_{n}$ is a subfield of $\Lambda _{n+1}$, and $\mathcal{G}(\Lambda
_{n}/E _{0})$ is cyclic of order $p ^{n}$. Using the fact that
$\Lambda _{n}$ is tamely and totally ramified over $E _{0}$, one
obtains from Lemma \ref{lemm3.1} and \cite{L2}, Ch. II,
Proposition~12, that $\Lambda _{n}/E _{0}$ possesses a primitive
element which is a $p ^{n}$-th root of some element $\pi _{n} \in E
_{0} ^{\ast }$, such that $v(\pi _{n}) > 0$ and $v(\pi _{n})$
generates the cyclic group $v(E _{0})$.
\par
In order to complete the proof of Theorem \ref{theo2.1} we use the
fact that $\Lambda _{n}/E _{0}$ is a Galois extension. Therefore,
the preceding observation shows that $\Lambda _{n}$ is a root field
over $E _{0}$ of the polynomial $X ^{p ^{n}} - \pi _{n}$. This
implies $\varepsilon _{n} \in \Lambda _{n}$ and the extension $E
_{0}(\varepsilon _{n})/E _{0}$ is totally ramified. However, it has
already been proved that $\varepsilon _{n} \in E _{0,{\rm ur}} = K
_{\rm ur}$, so it follows from (3.3) (a) that $\varepsilon _{n} \in
E _{0}$, for all $n \in \mathbb{N}$. Theorem \ref{theo2.1} is
proved.
\par
\medskip
\begin{coro}
\label{coro4.2} Let $(K, v)$ be an {\rm HDV}-field and $E/K$ a
Galois extension satisfying the conditions of Theorem \ref{theo2.1}.
Assume that $K$ is an almost perfect field, {\rm char}$(\widehat K)
= q$, and $v(E) \neq qv(E)$. Then $E$ is a $C _{1}$-field.
\end{coro}
\par
\medskip
\begin{proof}
Under the conditions of Theorem \ref{theo2.1}, $\widehat K$ is
finite and we have dim$(E) \le 1$, so it follows from the inequality
$v(E) \neq qv(E)$ that $\widehat E$ is algebraically closed. This
implies $K _{\rm ur}$ is $K$-isomorphic to a subfield of $E$. At the
same time, the assumption that $K$ is almost perfect allows us to
deduce from Lang's theorem that $K _{\rm ur}$ is a $C _{1}$-field.
Therefore, algebraic extensions of $K _{\rm ur}$ are also $C
_{1}$-fields, which proves our assertion.
\end{proof}
\par
\medskip
\begin{rema}
\label{rema4.2} It is known that $\mathbb{Q}$ has a unique
$\mathbb{Z} _{p}$-extension $\Gamma _{p}$ in $\mathbb{Q} _{\rm
sep}$, for each $p \in \mathbb{P}$ (and $\Gamma _{p}$ is included in
$\mathbb{Q}(\Theta _{p})$, where $\Theta _{p}$ is the set of all
roots of unity in $\mathbb{Q} _{\rm sep}$ of $p$-primary degrees).
Also, it is clear from Galois theory that the compositum $\Gamma $
of fields $\Gamma _{p}$, $p \in \mathbb{P}$, is a Galois extension
of $\mathbb{Q}$ with $\mathcal{G}(\Gamma /\mathbb{Q}) \cong \prod
_{p \in \mathbb{P}} \mathbb{Z} _{p}$. Identifying $\mathbb{Q} _{\rm
sep}$ with its isomorphic copy in $\mathbb{Q} _{q}$, for an
arbitrary fixed $q \in \mathbb{P}$, and using the decomposition law
in cyclotomic extensions of $\mathbb{Q}$ (see \cite{CF}, Ch. III,
Lemmas~1.3, 1.4), one obtains by the method of proving Theorem
\ref{theo2.1} that $\Gamma \cap \mathbb{Q} _{q} = \mathbb{Q}$,
$\mathbb{Q} _{q}\Gamma /\mathbb{Q} _{q}$ is a Galois extension with
$\mathcal{G}(\mathbb{Q} _{q}\Gamma /\mathbb Q _{q}) \cong
\mathcal{G}(\Gamma /\mathbb{Q})$, dim$(\mathbb{Q} _{q}\Gamma ) \le
1$ and $s(\mathbb{Q} _{q}\Gamma ) = \{q\}$. In addition, if $v _{q}$
is a valuation of $\Gamma $ extending the standard $q$-adic
valuation, say $\omega _{q}$, of $\mathbb{Q}$, then $v _{q}(\Gamma )
= \{r \in \mathbb{Q}\colon q ^{\rho } \in \omega _{q}(\mathbb{Q}),$
for some $\rho = \rho (r) \in \mathbb{N} \cup \{0\}\}$.
\end{rema}
\par
\medskip
\begin{prop}
\label{prop4.3}
Let $(K, v)$ be an {\rm HDV}-field with $\widehat K$ finite, and let
$T$ be the compositum of tamely ramified finite extensions of $K$ in
$K _{\rm ab}$. Then $T/K$ is a Galois extension with
$\mathcal{G}(T/K)$ isomorphic to $\mathcal{G}_{\widehat K} \times
\widehat K ^{\ast }$; hence, the valuation $v _{T}$ of $T$ extending
$v$ is discrete.
\end{prop}
\par
\medskip
\begin{proof}
Clearly, $K _{\rm ab}/K$ is a Galois extension with $\mathcal{G}(K
_{\rm ab}/K)$ which ensures that so is $L/K$ whenever $L \in I(K
_{\rm ab}/K)$. Specifically, $T/K$ is Galois, so it remains to be
proved that $\mathcal{G}(T/K)$ and $v _{T}$ have the properties
claimed by Proposition \ref{prop4.3}. It is well-known that
$\mathcal{G}_{\widehat K}$ is isomorphic to the topological group
product $\prod _{p \in \mathbb{P}} \mathbb{Z} _{p}$, whence, it is a
projective profinite group, in the sense of \cite{S1}. Since, by
(3.4) (b), $\mathcal{G}(K _{\rm ur}/K) \cong \mathcal{G}_{\widehat
K}$, this enables one to deduce from Galois theory that $T$ equals
the compositum $K _{\rm ur}T _{0}$, for some $T _{0} \in I(T/K)$,
such that $K _{\rm ur} \cap T _{0} = K$. Hence, by (3.3) (b) and
\cite{TW}, Proposition~A.17, finite extensions of $K$ in $T _{0}$
are tamely and totally ramified, which implies $\widehat T _{0} =
\widehat K$ and $v(T _{0})/v(K)$ is an abelian torsion group without
elements of order $p$. Taking now into account that $(K, v)$ is an
HDV-field, and using Lemma \ref{lemm3.1} and \cite{L2}, Ch. II,
Proposition~12, one obtains that, for each finite extension $K _{1}$
of $K$ in $T _{0}$, $K$ and $\widehat K$ contain primitive roots of
unity of degree $[K _{1}\colon K]$. Also, it becomes clear that $K
_{1}/K$ is Galois, $\mathcal{G}(K _{1}/K)$ is cyclic of order $n$,
and $\mathcal{G}(K _{1}/K) \cong v(K _{1})/v(K)$. As $\widehat K$ is
finite, this leads to the conclusion that $v(T _{0})/v(K)$ is cyclic
of order $q - 1$, and $T _{0}/K$ is a Galois extension with
$\mathcal{G}(T _{0}/K) \cong v(T _{0})/v(K)$, where $q$ is the
cardinality of $\widehat K$. Thus it follows that $[T\colon K _{\rm
ur}] = [T _{0}\colon K] = q - 1$, which ensures that $(T, v _{T})$
is an HDV-field (cf. \cite{E3}, Corollary~14.2.2). Note finally that
$T = K _{\rm ur}T _{0}$ is a Galois extension of $K$ with
$\mathcal{G}(T/K)$ isomorphic to the topological direct products
$\mathcal{G}(K _{\rm ur}/K) \times \mathcal{G}(T _{0}/K)$ and
$\mathcal{G}_{\widehat K} \times \widehat K ^{\ast }$
($\mathcal{G}(T _{0}/K)$ and $\widehat K ^{\ast }$ are viewed as
discrete topological group), so Proposition \ref{prop4.3} is proved.
\end{proof}
\par
\medskip
\begin{coro}
\label{coro4.4}
Under the hypotheses of Proposition \ref{prop4.3}, let $E$ be an
extension of $K$ in $T$. Then {\rm dim}$(E) \le 1$ if and only if $E
\in I(T/K _{\rm ur})$; when this holds and char$(K) > 0$, $E$ is a $C
_{1}$-field.
\end{coro}
\par
\medskip
\begin{proof}
Proposition \ref{prop4.3} and our assumptions show that $(E, v _{E})$
is an {\rm HDV}-field, so the former part of our conclusion follows from
the former one and Corollary \ref{coro4.2}. The noted property of
$(E, v _{E})$ also enables one to deduce the former part of Corollary
\ref{coro4.4} from Theorem \ref{theo2.1} and Lemma \ref{lemm3.1}.
\end{proof}
\par
\medskip
\begin{coro}
\label{coro4.5} Let $K$ be a finite extension of $\mathbb{Q}$ and $L
_{1}$ an extension of $K$ in $K _{\rm ab}$ of degree $p \in
\mathbb{P}$. Suppose that $L _{1}/K$ is of infinite height, i.e.
there are fields $L _{n} \in I(K _{\rm ab}/L _{1})$, $n \in
\mathbb{N}$, such that $L _{n}/K$ is a cyclic extension of degree $p
^{n}$, for each index $n$. Then one of the following conditions
holds:
\par
{\rm (a)} $L _{1}/K$ is wildly ramified; more precisely, if $L
_{1}/K$ is totally ramified relative to a valuation $v$ of $K$, then
$v(p) > 0$;
\par
{\rm (b)} $L _{1} \in I(H(K)/K)$, where $H(K)$ is the Hilbert class
field of $K$; in this case, $p$ divides the class number of $K$.
\end{coro}
\par
\medskip
\begin{proof}
Let $w$ be an arbitrary discrete valuation of $K$, such that $w(p) =
0$, i.e. the characteristic of the residue field of $(K, w)$ is
different from $p$. Identifying $K _{\rm ab}$ with its
$K$-isomorphic copy in $K _{w,{\rm sep}}$, put $L _{n} ^{\prime } =
L _{n}K _{w}$, for each $n \in \mathbb{P}$. It is easy to see that
if $L _{1}/K$ is totally ramified relative to $w$, then so is $L
_{1} ^{\prime }/K _{w}$ (relative to the continuous prolongation
$\bar w$ of $w$ on $K _{w}$). In view of Galois theory and
\cite{TW}, Proposition~A.17, these observations ensure that the
extensions $L _{n} ^{\prime }/K _{w}$, $n \in \mathbb{N}$, are both
cyclic of degree $p ^{n}$ and tamely and totally ramified. Therefore,
$\bar w(T)$ must include $\bar w(K _{w}) = w(K)$ as a subgroup of
infinite index, $T$ being the compositum of tame finite extensions of
$K _{w}$ in $K _{w,{\rm ab}}$. Since $w(K) \neq \{0\}$, and by
Proposition \ref{prop4.3}, $\bar w(T)$ is cyclic, this is a
contradiction due to the assumption that $L _{1}/K$ is totally
ramified relative to $w$, so Corollary \ref{coro4.5} (a) is proved.
\par
We turn to the proof of Corollary \ref{coro4.5} (b). The latter part
of our assertion follows from the former one and well-known general
properties of $H(K)$. Also, it is clearly sufficient to prove the
former part of Corollary \ref{coro4.5} (b), under the hypothesis that
$L _{1}/K$ is unramified relative to $v$, for any discrete
valuation $v$ of $K$ (and each valuation $v _{1}$ of $L _{1}$
extending $v$). When $p > 2$ or $K$ is a nonreal field, our assertion
follows from the definition of $H(K)$, so we consider only the case
where $p = 2$ and $K$ is formally real. Fix an Archimedean absolute
value $\omega $ of $K$ so that the completion $K _{\omega }$ of $K$
with respect to the topology induced by $\omega $ be isomorphic to
$\mathbb{R}$. It is easily obtained from Galois theory that if $L
_{1} \otimes _{K} K _{\omega }$ is a field, then so must be $L _{2}
\otimes _{K} K _{\omega }$. This, however,
means that $(L _{1} \otimes _{K} K _{\omega })/K _{\omega }$ must be
a quartic extension, which is impossible. Therefore, $L _{1} \otimes
_{K} K _{\omega }$ is not a field, which implies that if $\omega
_{1}$ is an absolute value of $L _{1}$ extending $\omega $, then the
completion of $L _{1}$ with respect to the topology of $\omega _{1}$
is isomorphic to $\mathbb{R}$ (see \cite{CF}, Ch. II, Theorem~10.2).
Hence, $L _{1}/K$ is unramified, and since $L _{1} \in I(K _{\rm
ab}/K)$, it follows that $L _{1} \in I(H(K)/K)$, so Corollary
\ref{coro4.5} is proved.
\end{proof}
\par
\medskip
\begin{prop}
\label{prop4.6} Let $(K, v)$ be an {\rm HDV}-field and $E/K$ a
Galois extension. Suppose that {\rm char}$(K) = 0$, $\widehat K$ is
finite, {\rm dim}$(E) \le 1$, and $\mathcal{G}(E/K)$ is abelian.
Then there exists $\Lambda  \in I(E/K)$ with $\mathcal{G}(\Lambda
/K)$ isomorphic to the topological group product $\prod _{p \in
\mathbb{P}} \mathbb{Z} _{p}$.
\end{prop}
\par
\smallskip
\begin{proof}
In view of Galois theory, it suffices to prove that $K$ has a
$\mathbb{Z} _{p}$-extension $\Lambda _{p}$ in $E$, for each $p \in
\mathbb{P}$. The existence of $\Lambda _{p}$, $p \in \mathbb{P}
\setminus \{q\}$, where $q = {\rm char}(\widehat K)$, is implied by
(3.4), the structure of $\mathcal{G}_{\widehat K}$, and Theorem
\ref{theo2.1}. It remains to be seen that $K$ has a $\mathbb{Z}
_{q}$-extension $\Lambda _{q}$ in $E$. Clearly $p$ does not divide
the degree of any finite extension of $E \cap K(p)$ in $E$. Hence,
by Lemma \ref{lemm3.3}, $[(E \cap K(p))\colon K] = \infty $.
Moreover, it follows from Krasner's lemma (cf. \cite{L2}, Ch. II,
Propositions~3, 4) and Galois theory that $X(E/K) _{q}$ is an
infinite group with finitely many elements of order $q ^{n}$, for
each $n \in \mathbb{N}$. Using Galois theory and applying the
following lemma, one proves the existence of $\Lambda _{q}$.
\end{proof}
\par
\smallskip
\begin{lemm}
\label{lemm4.7} Assume that $A$ is an infinite abelian torsion
$p$-group, for some $p \in \mathbb{P}$, such that the subgroup
$_{p}A = \{a \in A\colon pa = 0\}$ is finite. Then $A$ possesses a
subgroup isomorphic to $\mathbb{Z}(p ^{\infty })$.
\end{lemm}
\par
\smallskip
\begin{proof}
Arguing by induction on $n$, one obtains first that $A$ contains
finitely many elements of order $p ^{n}$, for each $n \in
\mathbb{N}$. Consider a sequence $\overline C = C _{n}$, $n \in
\mathbb{N}$, of cyclic subgroups of $A$, chosen so that $C _{n}$ be
of order $p ^{n}$, for each index $n$. It follows from the finitude
of $_{p}A$ that there is a subsequence $\overline C _{1} = C
_{1,n}$, $n \in \mathbb{N}$, of $\overline C$, such that the groups
$C _{1,n}$, $n \in \mathbb{N}$, share a common subgroup $H _{1}$ of
order $p$. Similarly, there exist subsequences $\overline C _{k}$,
$k \in \mathbb{N}$, of $\overline C$, such that $\overline C _{k} =
C _{k,n}$, $n \in \mathbb{N}$, the groups $C _{k,n}$, $n \in
\mathbb{N}$, have a common subgroup $H _{k}$ of order $p ^{k}$, and
$\overline C _{k+1}$ is a subsequence of $\overline C _{k}$, for
each $k$. Observing that $H _{k}$ is a subgroup of $H _{k+1}$, one
concludes that the union $H = \cup _{k=1} ^{\infty } H _{k}$ is a
subgroup of $A$ and there is a group isomorphism $H \cong \mathbb{Z}
(p ^{\infty })$.
\end{proof}
\par
\smallskip
Given an HDV-field $(K, v)$ with $\widehat K$ finite and a Galois extension
$E/K$, the conclusion of Proposition \ref{prop4.6} need not be true
without the assumptions that char$(K) = 0$ and $\mathcal{G}(E/K)$ is
abelian. For example, if char$(K) = q > 0$, then $K(q)$ includes
infinitely many subfields $\Theta _{n}$, $n \in \mathbb{N}$, that
are Artin-Schreier extensions of $K$ of degree $q$. This implies the
compositum $\Lambda _{q}$ of the fields $\Theta _{n}$, $n \in
\mathbb{N}$, is a Galois extension of $K$ with $\mathcal{G}(\Lambda
_{q}/E)$ infinite abelian of period $q$. Assuming that $\Lambda _{p}$
is the $\mathbb{Z} _{p}$-extension of $K$ in $K _{\rm ur}$, for each
$p \in \mathbb{P} \setminus \{q\}$, one obtains that the compositum
$\Lambda $ of the fields $\Lambda _{p}$, $p \in \mathbb{P}$, is a
Galois extension of $K$ with $\mathcal{G}(\Lambda /K) \cong \prod _{p
\in \mathbb{P}} \mathcal{G}(\Lambda _{p}/K)$. Hence, dim$(\Lambda )
\le 1$ and $\mathcal{G}(\Lambda /K)$ is abelian without a quotient
group isomorphic to $\mathbb{Z} _{q}$. Therefore, there is no
$\mathbb{Z} _{q}$-extension of $K$ in $\Lambda $, proving that the
condition on char$(K)$ in Proposition \ref{prop4.6} is essential.
Similarly, if $(K, v)$ is an HDV-field with char$(K) = 0$ and
char$(\widehat K) = q$, then the conclusion of Proposition
\ref{prop4.6} need not be true when $\mathcal{G}(E/K)$ is nonabelian.
Omitting the details, note that a counter-example $E$ can be chosen
so as to satisfy the following: (i) $E \cap K _{\rm ab} \subset K
_{\rm ur}$ and $K$ does not possess a degree $q$ extension in $E \cap
K _{\rm ab}$; (ii) $E/K$ is a Galois extension with $\mathcal{G}(E/(E
\cap K _{\rm ab}))$ infinite abelian of period $q$.
\par
\medskip
\section{\bf Characterization of fields of dimension $\le 1$ within
the class of abelian and tame extensions of a global field}
\par
\medskip
Our main goal in this Section is to characterize the fields
pointed out in its title and thereby to prove Theorem \ref{theo2.2}
(a). First, we present some basic properties of maximal abelian tame
extensions of global fields, as follows:
\par
\medskip
\begin{prop}
\label{prop5.1}
Let $K$ be a global field and $\Theta $ the compositum of tame finite
extensions of $K$ in $K _{\rm ab}$. Then:
\par
{\rm (a)} All nontrivial valuations of $\Theta $ are discrete with
algebraically closed residue fields; in particular, {\rm
dim}$(\Theta ) \le 1$;
\par
{\rm (b)} For each nontrivial valuation $v$ of $K$, the compositum
$T$ of tamely ramified finite extensions of $K _{v}$ in $K _{v,{\rm
ab}}$ contains as a subfield a $K$-isomorphic copy $\Theta ^{\prime
}$ of $\Theta $; in addition, $\Theta ^{\prime }K _{v} = T$.
\end{prop}
\par
\medskip
\begin{proof}
Identifying $K _{\rm sep}$ with its $K$-isomorphic copy in $K
_{v,{\rm sep}}$, one may put $\Theta ^{\prime } = \Theta $. Using
also Galois theory and the tameness of $\Theta /K$, one obtains
further that $K _{\rm ab}K _{v} \subseteq K _{v,{\rm ab}}$ and
$\Theta .K _{v} \subseteq T$. Conversely, it follows from
Grunwald-Wang's theorem (see \cite{LR}) that each cyclic extension
$\Lambda $ of $K _{v}$ in $K _{v,{\rm ab}}$ equals $\Lambda _{0}K
_{v}$, for some cyclic extension $\Lambda _{0}$ of $K$ in $K _{\rm
ab}$ satisfying the divisibility conditions $[\Lambda \colon K _{v}] \mid
[\Lambda _{0}\colon K] \mid 2[\Lambda \colon K _{v}]$. Moreover,
Grunwald-Wang's theorem guarantees that if $\Lambda \in I(T/K)$,
then $\Lambda _{0}$ can be chosen so as to lie in $I(\Theta /K)$.
This implies the inclusions $T \subseteq \Theta K _{v}$ and $K
_{v,{\rm ab}} \subseteq K _{\rm ab}K _{v}$. These observations prove
Proposition \ref{prop5.1} (b). Consider now the fields $\Theta $ and
$\Theta _{1} = \Theta .K _{1}$, where $K _{1}$ is the maximal
separable (algebraic) extension of $K$ in $K _{v}$. Denote by
$\theta $ and $\theta _{1}$ the valuations of $\Theta $ and $\Theta
_{1}$, respectively, induced by the valuation of $T = \Theta K _{v}$
extending the continuous prolongation $\bar v$ of $v$ on $K _{v}$.
Clearly, $\theta $ and $\theta _{1}$ are discrete ($\bar v _{T}$ is
discrete, by Proposition \ref{prop4.3}, and $\{0\} \neq v(K) \le
\theta (\Theta ) \le \theta _{1}(\Theta _{1}) \le \bar v(T)$), and
$(K _{1}, v _{1})$ is a Henselization of $(k, v)$, where $v _{1}$ is
the valuation of $K _{1}$ induced by $\bar v$. The latter ensures
that $(\Theta _{1}, \theta _{1})$ is a Henselization of $(\Theta ,
\theta )$ (cf. \cite{TW}, Proposition~A.30). Taking also into account
that $\Theta _{1}K _{v} = \Theta K _{v} = T$, one deduces from Lemma
\ref{lemm3.2} and Proposition \ref{prop4.3} that $\widehat \Theta
_{1} = \widehat \Theta \cong \widehat K _{\rm sep}$. In view of the
normality of $\Theta /K$, every valuation of $\Theta $ possesses the
obtained properties of $\theta $, so the former part of Proposition
\ref{prop5.1} (a) is proved. Note finally that $\Theta $ is a nonreal
field, i.e. $-1$ is presentable as a sum of the squares of finitely
many elements of $\Theta $. This is obvious in case char$(K) > 0$,
and if char$(K) = 0$, then $\Theta $ contains a square root
$\sqrt{-p}$,  for any $p \in \mathbb{P}$ satisfying $p \equiv 3 ({\rm
mod} \ 4)$ (whence, $-1$ is presentable as a sum of squares of four
elements of $\Theta $). Since, Corollary \ref{coro4.4} and the former
part of Proposition \ref{prop5.1} (a) yield Br$(T) = \{0\}$, for each
nontrivial valuation $v$ of $\Theta $, this allows to deduce
from Lemma \ref{lemm3.4} that Br$(\Theta ) = \{0\}$. Hence, by
\cite{FS}, Theorem~4, dim$(\Theta ) \le 1$, so Proposition
\ref{prop5.1} is proved.
\end{proof}
\par
\medskip
Theorem \ref{theo2.2} (a) is contained in Proposition \ref{prop5.1}
(a) and the following:
\par
\medskip
\begin{coro}
\label{coro5.2}
Let $K$ be a global field and $E$ a tame extension of $K$ in $K _{\rm
ab}$. Then all nontrivial Krull valuations of $E$ are discrete. In
addition, {\rm dim}$(E) \le 1$ if and only if $E$ is a nonreal field
and the residue fields of its discrete valuations are algebraically
closed.
\end{coro}
\par
\medskip
\begin{proof}
Let $\Theta $ be the compositum of tame finite extensions of $K$ in
$K _{\rm ab}$. Then $E \in I(\Theta /K)$ and every nontrivial
valuation $w$ of $E$ extends to a valuation $w'$ of $\Theta $,
so the former conclusion of Corollary \ref{coro5.2} follows from
Proposition \ref{prop5.1} (a). Note further that if $E$ is a formally
real field, i.e. $-1$ is not presentable as a sum of the squares of
finitely many elements of $E$, then the Hamiltonian quaternion
$E$-algebra is a division one, whence, Br$(E)$ contains an element of
order $2$. When $E$ is nonreal, the latter assertion of the corollary
can be deduced following the concluding part of the proof of
Proposition \ref{prop5.1}.
\end{proof}
\par
\medskip
The following result shows that if $K$ is a number field with Cl$(K)
= 1$, then $X(\Delta /K)$ decomposes into a direct sum of finite
cyclic groups, for each tame extension $\Delta $ of $K$ in $K _{\rm
ab}$. This need not be valid without the tameness condition; one may
take as a counter-example any nonreal field $\Delta ^{\prime } \in
I(K _{\rm ab}/\Gamma K)$ with $[\Delta ^{\prime }\colon \Gamma K] <
\infty $, where $\Gamma $ is defined in Remark \ref{rema4.2}.
\par
\smallskip
\begin{prop}
\label{prop5.3} Let $K$ be a number field, {\rm Cl}$(K)$ its class
number, $\Delta $ a tame extension of $K$ in $K _{\rm ab}$, and
$X(\Delta /K) _{p}$ the $p$-component of $X(\Delta /K)$, for each $p
\in \mathbb{P}$. Then:
\par
{\rm (a)} $X(\Delta /K)$ is a reduced abelian torsion group with
finitely many elements of infinite height;
\par
{\rm (b)} $X(\Delta /K) _{p}$ contains infinitely many elements of
order $p$ unless it is a finite group; $X(\Delta /K) _{p}$ is
infinite in case {\rm Br}$(\Delta ) _{p} = \{0\}$;
\par
{\rm (c)} $X(\Delta /K)$ decomposes into a direct sum of cyclic
$p$-groups, for every $p \in \mathbb{P}$ not dividing {\rm Cl}$(K)$.
\end{prop}
\par
\medskip
\begin{proof}
Clearly, $X(\Delta /K)$ is a countable abelian torsion group, and
Corollary \ref{coro4.5} shows that $X(\Delta /K) _{p}$ does not
contain nonzero elements of infinite height, for any $p \in
\mathbb{P}$, $p \nmid {\rm Cl}(K)$, so Proposition \ref{prop5.3} (c)
follows from Pr\"{u}fer's theorem (see \cite{F}, Theorem~5.3) and
Proposition \ref{prop5.3} (a). For the proof of former part of
Proposition \ref{prop5.3} (b), one may clearly assume that $X(\Delta
/K) _{p}$ is infinite. Suppose for a moment that $X(\Delta /K) _{p}$
contains only finitely many elements of order $p$. Then, by Lemma
\ref{lemm4.7}, $X(\Delta /K) _{p}$ must have a subgroup isomorphic to
$\mathbb{Z}(p ^{\infty })$. In view of Galois theory, this requires
the existence of a $\mathbb{Z} _{p}$-extension $\Lambda _{p}$ of $K$
in $\Delta $. Clearly, $\Lambda _{p}/K$ must preserve the tameness of
$\Delta /K$, whence, using repeatedly Corollary \ref{coro4.5}, one
concludes that $\Lambda _{p}$ must be a subfield of $H(K)$. On the
other hand, by class field theory (Furtw\"{a}ngler's theorem, see
\cite{CF}, Chs. IX and XI), $H(K)/K$ is a finite extension with
$[H(K)\colon K] = {\rm Cl}(K)$. The obtained contradiction proves the
former part of Proposition \ref{prop5.3} (b).
\par
For the rest of the proof of Proposition \ref{prop5.3} (b), note that
$X(\Delta /K) _{p}$ is finite if and only if $\Delta \cap K(p)$ is a
finite extension of $K$ (this is a well-known consequence of Galois
theory and Pontrjagin's duality, see  ). Since, by class field theory
(cf. \cite{We}, Ch. XIII, Sects. 3 and 6), Br$(Y) _{p'} \neq \{0\}$,
$p' \in \mathbb{P}$, for every global field $Y$, this implies the
latter part of Proposition \ref{prop5.3} (b).
\par
It remains for us to prove Proposition \ref{prop5.3} (a). Clearly,
the set $X _{0}(\Delta /K)$ of elements of $X(\Delta /K)$ of
infinite height forms a subgroup of $X(\Delta /K)$, so it follows from
Proposition \ref{prop5.3} (c) that it suffices to show that if Cl$(K)
> 1$, then the $p$-component $X _{0}(\Delta /K) _{p}$ of $X
_{0}(\Delta /K)$ is a finite group, for an arbitrary $p \in
\mathbb{P}$ dividing Cl$(K)$. Since $K$ has finitely many
valuations, up-to equivalence, with residue fields of characteristic
$p$, Corollary \ref{coro4.5} (a) and \cite{Ko}, Theorem~1.48 (see
also \cite{Wash}, Ch. 13) imply the existence of finitely many
degree $p$ extensions of $K$ in $K _{\rm ab}$ of infinite height.
In view of Galois theory, this means that $X _{0}(\Delta /K) _{p}$
contains finitely many elements of order $p$. Using Galois theory,
one also concludes there is $I _{p} \in I(\Delta /K)$ with $X(I
_{p}/K) \cong X _{0}(\Delta /K) _{p}$. This allows to obtain by the
method of proving Proposition \ref{prop5.3} (b) that $\mathbb{Z}(p
^{\infty })$ does not embed as a subgroup of $X _{0}(\Delta /K)
_{p}$. Hence, by Lemma \ref{lemm4.7}, $X _{0}(\Delta /K)  _{p}$ is
finite, so Proposition \ref{prop5.3} is proved.
\end{proof}
\par
\smallskip
\begin{rema}
\label{rema5.4} Let $K$ be a number field. Then Corollary
\ref{coro4.5} and \cite{Ko}, Theorem~1.48, imply that for any $p \in
\mathbb{P}$, $K$ has finitely many extensions in $K _{\rm ab}$ of
degree $p$ and infinite height. Hence, by Galois theory and the
proof of Lemma \ref{lemm4.7}, $X(K _{\rm ab}/K)$ contains finitely
many elements of order $p ^{n}$ and infinite height, for each $n \in
\mathbb{N}$. In view of the structure and injectivity of divisible
abelian torsion groups (cf. \cite{F}, Theorems~2.6, 3.1), this
recovers the proof of \cite{FS2}, Theorem~1, in the case of
a number ground field. It is not known whether the quotient group of 
$X(K _{\rm ab}/K)$ by its maximal divisible subgroup $DX(K _{\rm 
ab}/K)$ contains only finitely many elements of infinite height. By 
\cite{FS2}, Theorem~14, this holds for global function fields $K 
^{\prime }$.
\end{rema}
\par
\medskip
The maximal tame extension of $\mathbb{Q}$ in $\mathbb{Q} _{\rm ab}$ 
is described explicitly as follows:
\par
\medskip
\begin{prop}
\label{prop5.7} The maximal tame extension $\Theta $ of
$\mathbb{Q}$ in $\mathbb{Q} _{\rm ab}$ is generated over
$\mathbb{Q}$ by the set $\Psi $ of primitive $p$-th roots of unity
$\varepsilon _{p} \in \mathbb{Q} _{\rm sep}$, $p \in \mathbb{P}$.
\end{prop}
\par
\smallskip
\begin{proof}
Denote by $E$ the extension of $\mathbb{Q}$ generated by $\Psi $, and
let $\Gamma $ be the field defined in Remark \ref{rema4.2}. It is
well-known (see \cite{CF}, Ch. II, Theorem~3.2) that the natural
$p$-adic valuations $v _{p}$, $p \in \mathbb{P}$, form a system of
representatives of the equivalence classes of nontrivial Krull
valuations of $\mathbb{Q}$. Assume now that $\Phi $ is the set of
square-free odd integers $\ge 3$, and for each $n \in \Phi $, put
$\Lambda _{n} = \mathbb{Q}(\varepsilon _{n})$, where $\varepsilon
_{n} \in \mathbb{Q} _{\rm sep}$  is a primitive $n$-th root of unity;
also, let $\lambda _{n,p}$ be a valuation of $\Lambda _{n}$ extending
$v _{p}$, for each $p \in \mathbb{P}$. Then $\Lambda _{n}/\mathbb{Q}$
is a Galois extension with $\mathcal{G}(\Lambda _{n}/\mathbb{Q})
\cong (\mathbb{Z}/n\mathbb{Z}) ^{\ast }$ (cf. \cite{CF}, Ch. III,
Lemma~1.1). As $n \in \Phi $, this means that $\mathcal{G}(\Lambda
_{n}/\mathbb{Q})$ is isomorphic to the direct product $\prod _{p'}
\mathbb{F} _{p'} ^{\ast }$, indexed by the set of prime divisors of
$n$. It is therefore clear that, for any pair $\nu _{1}$, $\nu _{2} \in
\Phi $, $\Lambda _{\nu _{1}}\Lambda _{\nu _{2}} = \Lambda _{\nu }$,
where $\nu = {\rm lcm}[\nu _{1}, \nu _{2}]$; in particular $E$ equals
the union of $\Lambda _{n}$, $n \in \Phi $. To prove the inclusion $E
\subseteq \Theta $, it is sufficient to show that $\Lambda
_{n}/\mathbb{Q}$ is tamely ramified relative to $\lambda _{n,p}/v
_{p}$, for arbitrary fixed $n \in \Phi $, $p \in \mathbb{P}$. This is
implied by the following facts: $\lambda _{n,p}(\Lambda _{n}) = v
_{p}(\mathbb{Q})$, provided that $p \nmid n$; $v _{p}(\mathbb{Q})$ is
a subgroup of $\lambda _{n,p}(\Lambda _{n})$ of index $p - 1$ in case
$p \mid n$ (cf. \cite{CF}, Ch. III, Lemmas~1.3 and 1.4). Thus the
assertion that $E \subseteq \Theta $ becomes obvious as well as the
fact that every valuation $w _{p}$ of $E$ extending $v _{p}$ is
discrete with $w _{p}(E)$ including $v _{p}(\mathbb{Q})$ as a
subgroup of index $p - 1$.
\par
It remains to prove that $\Theta = E$. It follows from the
Kronecker-Weber theorem, Galois theory and well-known basic
properties of cyclotomic extensions of $\mathbb{Q}$ that $\mathbb{Q}
_{\rm ab} = E\Gamma (\sqrt{-1})$ and $E \cap \Gamma(\sqrt{-1}) =
\mathbb{Q}$. This implies $\mathcal{G}(\mathbb{Q} _{\rm ab}/E) \cong
\mathcal{G}(\Gamma (\sqrt{-1})/\mathbb{Q})$ and the map of
$I(\Gamma (\sqrt{-1})/\mathbb{Q})$ into $I(\mathbb{Q} _{\rm ab}/E)$,
by the rule $Y \to YE$, is bijective. Also, it is clear that
$[Y _{1}\colon \mathbb{Q}] = [Y _{1}E\colon E]$, for any finite
extension $Y _{1}$ of $\mathbb{Q}$ in $\Gamma (\sqrt{-1})$. We show
that if $Y _{1} \neq \mathbb{Q}$, then $Y _{1}/\mathbb{Q}$
is not tame. Indeed, for any divisor $p \in \mathbb{P}$ of
$[Y _{1}\colon \mathbb{Q}]$, $Y _{1}/\mathbb{Q}$ is wildly ramified
relative to $y _{p}/v _{p}$, where $y _{p}$ is a valuation of $Y
_{1}$ extending $v _{p}$ (apply \cite{CF}, Ch. II, Lemmas~1.3, 1.4).
Since the set of tame extensions of $\mathbb{Q}$ in $\mathbb{Q}
_{\rm ab}$ is closed under taking intermediate fields (cf.
\cite{L2}, Ch. II, Propositions~8, 13, and \cite{CF}, Theorem~10.2),
this proves that $E = \Theta $, as required.
\end{proof}
\par
\smallskip
Let us note that Proposition \ref{prop5.7} allows us to give an
alternative proof of Proposition \ref{prop5.1} (b). Indeed,
Theorem~4 of \cite{FS}, reduces our considerations to the special
case where $K = \mathbb{Q}$. Then our assertion can be obtained by
applying the following statement (which in turn is implied by the
Lemma in \cite{Chat}, page~131): for any fixed pair $(l, n) \in
\mathbb{P} \times \mathbb{N}$, there are infinitely many $q _{n} \in
\mathbb{P} \setminus \{p\}$, such that the extension $\mathbb{F}
_{p}(\gamma _{n})/\mathbb{F} _{p}$, where $\gamma _{n}$ is a
primitive $q _{n}$-th root of unity in $\mathbb{F} _{p,{\rm sep}}$,
has degree $[\mathbb{F} _{p}(\gamma _{n})\colon \mathbb{F} _{p}]$
divisible by $l ^{n}$.
\par
\medskip
\begin{coro}
\label{coro5.8} Assume that $\Theta _{\Pi }$ is an extension of
$\mathbb{Q}$ obtained by adjunction of primitive $p$-th roots of
unity $\varepsilon _{p} \in \mathbb{Q} _{\rm sep}$, where $p$ runs
across the complement $\mathbb{P} \setminus \Pi $, for an arbitrary
finite subset $\Pi $ of $\mathbb{P}$. Then {\rm dim}$(\Theta _{\Pi
}) \le 1$.
\end{coro}
\par
\smallskip
\begin{proof}
Let $\Theta $ be the field defined in Proposition \ref{prop5.7}.
Then $\Theta /\Theta _{\Pi }$ is a finite extension. More precisely,
if $\Pi = \{p _{1}, \dots , p _{s}\}$, then $[\Theta \colon \Theta
_{\Pi }] = \prod _{i=1} ^{s} (p _{i} - 1)$. This means that
$\mathcal{G}_{\Theta }$ is an open subgroup of $\mathcal{G}_{\Theta
_{\Pi }}$ of index $\prod _{i=1} ^{s} (p _{i} - 1)$. Note also that
$\Theta _{\Pi }$ is a nonreal field, because an ordered field does
not contain a primitive root of unity of any degree greater than
$2$. It is therefore clear from Galois theory that
$\mathcal{G}_{\Theta _{\Pi }}$ is a torsion-free group. Consider now
an arbitrary nontrivial valuation $w$ of $\Theta _{\Pi }$, fix a
valuation $w'$ of $\Theta $ extending $w$, and denote by $\widehat
\Theta _{\Pi }$ and $\widehat \Theta $ the residue fields of
$(\Theta _{\Pi }, w)$ and $(E, w')$, respectively. Then $w$ is
discrete and $\widehat \Theta /\widehat \Theta _{\Pi }$ is a finite
extension. Moreover, it follows from Propositions \ref{prop5.1} (a)
and \ref{prop5.7} that $\widehat \Theta $ is algebraically closed.
In view of Galois theory, this yields $\widehat \Theta = \widehat
\Theta _{\Pi }$. It is now easy to see that dim$(\Theta _{\Pi }) \le
1$.
\end{proof}
\par
\medskip
In the setting of Propositions \ref{prop5.1} (a) and \ref{prop5.7},
the question of whether a field $\Delta \in I(\Theta /\mathbb{Q})$
with dim$(\Delta ) \le 1$ is of type $C _{1}$ is widely open.
Proposition \ref{prop5.1} and \cite{L1}, Theorem~10, show that the
completion $\Delta _{v}$ is a $C _{1}$-field, for any nontrivial
valuation $v$ of $\Delta $. Thus our question concerning $\Delta $
is equivalent to the problem of finding whether the Hasse principle
applies to $\Delta $-forms in more variables than their degrees. It
is worth mentioning that, for each $q \in \mathbb{P}$, $\Theta $ can
be viewed as an analog in characteristic $q$ to the rational
function field $\mathbb{F} _{q,{\rm sep}}(X)$ in one variable over
$\mathbb{F} _{q,{\rm sep}}$ (which is of type $C _{1}$, by Tsen's
theorem). In contrast to the extension $\Theta /\mathbb{Q}$,
however, $\mathbb{F} _{q,{\rm sep}}(X)/\mathbb{F} _{q}(X)$ has no
intermediate field $\Delta _{q} \neq \mathbb{F} _{q,{\rm sep}}(X)$
with dim$(\Delta _{q}) \le 1$, for any $q$.
\par
\medskip
\section{\bf Proof of Theorem \ref{theo2.2}}
\par
\medskip
Assume that $E$ and $K$ are Galois extensions of $\mathbb{Q}$ in
$\mathbb{Q} _{\rm sep}$, such that $[K\colon \mathbb{Q}] < \infty $,
$E \in I(K _{\rm ab}/K)$, dim$(E) \le 1$, and $E/K$ is tame. We prove
the existence of an $n$-variate $E$-form $N _{n}$ of degree $n$,
violating the Hasse principle, by showing that there is an abelian 
unramified Galois extension $E _{n}$ of $E$ in $\mathbb{Q} _{\rm 
sep}$ with $[E _{n}\colon E] = n$. Then one may take as $N _{n}$ the 
norm form of $E _{n}/E$ associated with any $E$-basis of $E _{n}$ 
(since $N _{n}$ has no nontrivial zero over $E$, whereas Propositions 
\ref{prop5.1} (a) and \ref{prop5.7} imply $E _{n}/E$ is split 
relative to $v$, whence $N _{n}$ decomposes over $E _{v}$ into a 
product of $n$ linear forms in $n$ variables and so it has a 
nontrivial zero over $E _{v}$). The existence of the extension $E 
_{n}/E$ is ensured by the following result.
\par
\medskip
\begin{lemm}
\label{lemm6.1}
Let $K$ be a finite extension of $\mathbb{Q}$ in $\mathbb{Q} _{\rm 
ab}$, $H(K)$ the Hilbert class field of $K$, and $p \in \mathbb{P}$ a 
divisor of $[(H(K) \cap \mathbb{Q} _{\rm ab})\colon \mathbb{Q}$. Then 
$p \mid [K\colon \mathbb{Q}]$.
\end{lemm}
\par
\medskip
\begin{proof}
It follows  from the normality of $K/\mathbb{Q}$ and the definition 
of $H(K)$ that $H(K)/\mathbb{Q}$ is a finite Galois extension. Note 
further that our assumptions ensure the existence of a Galois 
extension $L _{p}$ of $\mathbb{Q}$ in $H(K)$ of degree $[L _{p}\colon 
\mathbb{Q}] = p$; this is clear from Galois theory and the general 
structure of finite abelian groups. Hence, by the Hermite-Minkowski 
theorem (cf. \cite{L2}, Ch. V, Sect. 4), there exists $\lambda \in 
\mathbb{P}$, such that $L _{p}/\mathbb{Q}$ is totally ramified at the 
$\lambda $-adic valuation $\omega _{\lambda }$ of $\mathbb{Q}$. This 
implies $\lambda \mid e(H(K)/\mathbb{Q}) _{\omega _{\lambda }}$, 
which indicates that $p \mid [K\colon \mathbb{Q}]$ (since $H(K)/K$ is 
unramified). 
\end{proof}
\par
\medskip
\begin{lemm}
\label{lemm6.2}
Let $K _{1}$ and $K _{2}$ be finite Galois extensions of $\mathbb{Q}$ 
in $\mathbb{Q} _{\rm sep}$ of $p$-primary degrees, for a given $p \in 
\mathbb{P}$, and $H(K _{2}) \in I(\mathbb{Q} _{\rm sep}/\mathbb{Q})$, be 
the Hilbert class fields of $K _{1}$ and $K _{2}$, respectively. 
Suppose that $K _{1} \cap K _{2} = \mathbb{Q}$, and for $j = 1, 2$, 
denote by $H _{j}$ the maximal subfield of $H(K _{j})$ with respect 
to the property that $[H _{j}\colon K _{j}]$ is not divisible by $p$. 
Then $H _{1}$, $H _{2}$ and $H _{1}H _{2}$ are Galois extensions of 
$\mathbb{Q}$, such that $H _{1} \cap H _{2} = \mathbb{Q}$, $H _{1}H 
_{2} \cap \mathbb{Q} _{\rm ab} = K _{1}K _{2}$, and $\mathcal{G}(H 
_{1}H _{2}/\mathbb{Q})$ is isomorphic to the direct product 
$\mathcal{G}(H _{1}/\mathbb{Q}) \times \mathcal{G}(H 
_{2}/\mathbb{Q})$.
\end{lemm}
\par
\medskip
\begin{proof}
Denote for brevity by $H _{0}$ and $H _{3}$ the fields $H _{1} \cap H 
_{2}$ and $H _{1}H _{2}$, respectively, and put $K _{3} = K _{1}K _{2}$.
Our assumptions, the normality of $K _{1}$ and $K _{2}$ over 
$\mathbb{Q}$, and properties of $H(K _{1})/K _{1}$, $H(K _{2})/K 
_{2}$ available by definition, ensure that $H _{1}/\mathbb{Q}$ and $H 
_{2}/\mathbb{Q}$ are Galois extensions, so it follows from Galois theory 
that $H _{3}/\mathbb{Q}$ is Galois as well. They also show that 
$\mathcal{G}(H _{j}/K _{j})$ is abelian and $\mathcal{G}(H 
_{j}/\mathbb{Q})$, $\mathcal{G}(H _{2}/\mathbb{Q})$ is 
metabelian, for $j = 1, 2, 3$. Taking further into account $p \nmid 
[H _{j}\colon K _{j}]$, $j = 1, 2, 3$, one concludes that, for each 
$j$, $K _{j}$ contains as a subfield every extension of $\mathbb{Q}$ 
in $H _{j}$ of $p$-primary degree. At the same time, it is clear from 
Galois theory that $H _{0}/\mathbb{Q}$ is a Galois extension and 
$\mathcal{G}(H _{0}/\mathbb{Q})$ is a homomorphic image of 
$\mathcal{G}(H _{j}/\mathbb{Q})$, $j = 1, 2$; in particular, 
$\mathcal{G}(H _{0}/\mathbb{Q})$ is metabelian. These facts and the 
equality $K _{1} \cap K _{2} = \mathbb{Q}$ indicate that if 
$H _{0} \neq \mathbb{Q}$, then $\mathbb{Q}$ has a Galois 
extension $K _{0}$ in $H _{0}$ of prime degree $[K _{0}\colon 
\mathbb{Q}] = l \neq p$. This, however, contradicts Lemma 
\ref{lemm6.1} and thereby proves that $H _{0} = \mathbb{Q}$. Now it 
follows from Galois theory that $\mathcal{G}(H 
_{3}/\mathbb{Q}) \cong \mathcal{G}(H _{1}/\mathbb{Q}) \times 
\mathcal{G}(H _{2}/\mathbb{Q})$, which implies $\mathcal{G}(H _{3}/K 
_{3})$ equals the commutator subgroup of $\mathcal{G}(H 
_{3}/\mathbb{Q})$ and $K _{3} = H _{3} \cap \mathbb{Q} _{\rm ab}$. 
Lemma \ref{lemm6.2} is proved.
\end{proof}
\par
\medskip
We are now in a position to prove (2.5) (and Theorem \ref{theo2.2} 
(b)), for a given integer $n \ge 2$. In view of the general structure 
of finite abelian groups, one may consider only the special case 
where $n$ is odd or $n = 2 ^{\mu }$, for some $\mu \in \mathbb{N}$. 
Our argument relies on the existence of Galois extensions $Q _{\nu }, 
Y _{\nu }$, $\nu \in \mathbb{N}$, with class numbers Cl$(Q _{\nu })$ 
and Cl$(Y _{\nu })$ divisible by $n$, and with $[Q _{\nu }\colon 
\mathbb{Q}] = 2$ and $[Y _{\nu }\colon \mathbb{Q}] = 3$, for every 
index $\nu $ (see \cite{Yam} and \cite{Uch}, respectively). Clearly, 
the sequence $Q _{\nu }, Y _{\nu }$, $m \in \mathbb{N}$, can be 
chosen so that $[Q _{\nu _{1}} \dots Q _{\nu _{k}}\colon \mathbb{Q}] 
= 2 ^{k}$ and $[Y _{\nu _{1}} \dots Y _{\nu _{k}}\colon \mathbb{Q}] = 
3 ^{k}$ whenever $k \ge 2$ is an integer and $\nu _{1} < \dots < \nu 
_{k}$. For each $\nu \in \mathbb{N}$, denote by $C _{\nu }$ the field 
$Q _{\nu }$ if $2 \nmid n$, and put $C _{\nu } = Y _{\nu }$ when $n = 
2 ^{\mu }$; also, let $H _{\nu }$ be an extension of $C _{\nu }$ in 
the Hilbert class field $H(C _{\nu }) \in I(\mathbb{Q}_{\rm sep}/C 
_{\nu })$ of degree $[H _{\nu }\colon C _{\nu }]$ not divisible by 
$[C _{\nu }\colon \mathbb{Q}]$. Proceeding by induction on $k \in 
\mathbb{N} \setminus \{1\}$, and using Galois theory, Lemma 
\ref{lemm6.2} and properties of Hilbert class fields available by 
definition, one obtains that $H _{\nu _{1}} \dots H _{\nu _{k}} \cap 
\mathbb{Q} _{\rm ab} = C _{\nu _{1}} \dots C _{\nu _{k}}$, $[H _{\nu 
_{1}} \dots H _{\nu _{k}}\colon \mathbb{Q}] = \prod _{j=1} ^{k} [H 
_{\nu _{j}}\colon \mathbb{Q}]$ and $[H _{\nu _{1}} \dots H _{\nu 
_{k}}\colon C _{\nu _{1}} \dots C _{\nu _{k}}] = \prod _{j=1} ^{k} [H 
_{\nu _{j}}\colon C _{\nu _{j}}]$, provided that $k \ge 2$ and $\nu 
_{1} < \dots < \nu _{k}$. It is now easy to see that if $H _{\nu }$, 
$\nu \in \mathbb{N}$, are chosen so that $[H _{\nu }\colon C _{\nu }] 
= n$, for every $\nu $ (which is possible, since $n \mid {\rm Cl}(C 
_{\nu })$), then the fields $E _{\nu } = H _{\nu }E$, $\nu \in 
\mathbb{N}$, have the properties required by (2.5), so the $E$-norms 
$N _{\nu }$ of $E _{\nu }/E$, $\nu \in \mathbb{N}$, are pairwise 
nonequivalent. 
\par
\medskip
\begin{coro}
\label{coro6.3} Assume that $K$ and $E$ satisfy the conditions
of Theorem 2.2  (b), and let $\Delta $ be an extension of $E$ in $K 
_{\rm ab}$. Then there exist abelian extensions $\Delta _{\nu }$, 
$\nu \in \mathbb{N}$, of $\Delta $ in $\mathbb{Q} _{\rm sep}$ with 
the following properties:
\par
{\rm (a)} $[\Delta _{\nu }\colon \Delta ] = n$, for every index $\nu $;
\par
{\rm (b)} The tensor product of the $\Delta $-algebras $\Delta _{\nu 
}$, $\nu \in \mathbb{N}$, is a field;
\par
{\rm (c)} If $N _{\nu }$ is a norm form of $\Delta _{\nu }/\Delta $, 
for each $\nu $, then $N _{\nu }$, $\nu \in \mathbb{N}$, are pairwise
nonequivalent $\Delta $-forms of degree $n$, violating the Hasse 
principle.
\end{coro}
\par
\medskip
\begin{proof} 
Clearly, one may take $E _{\nu }$, $\nu \in \mathbb{N}$, as in the 
proof of Theorem \ref{theo2.2} (b) and put $\Delta _{\nu } = \Delta E 
_{\nu }$, for each index $\nu $.
\end{proof}
\par
\medskip
Corollary \ref{coro6.3}, applied to the case of $K = \mathbb{Q}$,
shows that the latter part of the conclusion of Theorem \ref{theo2.2}
(b), remains valid, for any number field $K$, when the assumption
that $E = \Theta $ is replaced by the one that $E = \Theta _{\Pi }K$,
where $\Pi $ is a finite subset of $\mathbb{P}$ and $\Theta _{\Pi }$
is defined as in Corollary \ref{coro5.8}.
\par
\medskip
\emph{Acknowledgements.} This research has partially been supported
by Grant KP-06 N 32/1 of 07.12. 2019 of the Bulgarian National
Science Fund.
\par

\end{document}